\documentclass[12pt]{article}

\usepackage{a4wide}
\usepackage{amssymb}
\usepackage{amsfonts}
\usepackage{amsmath}
\input xy
\xyoption{arrow} \xyoption{matrix}

\date{}

\newtheorem{proposition}{Proposition}[section]
\newtheorem{theorem}[proposition]{Theorem}
\newtheorem{lemma}[proposition]{Lemma}

\newtheorem{corollary}[proposition]{Corollary}

\def\der{\partial }

\def\nFM0{{\nu }_{F,M_0}}
\def\nFN0{{\nu }_{F,N_0}}
\def\nGN0{{\nu }_{G,N_0}}

\def\N0{ {\bf N}_0 }

\def\g{\gamma}

\def\ra{\rightarrow}

\def\Xpm{X^{\pm }}

\def\s{\sigma}
\def\Z{\mathbb{Z}}

\def\l1{{\lambda}_1}

\def\a{\alpha}
\def\a0{ {\alpha }_0}
\def\a1{ {\alpha }_1}

\def\l{\lambda}
\def\o{\omega}

\def\nFGM0{{\nu }_{F,G,M_0}}

\def\nFN0{{\nu}_{F,N_0}}

\def\sm{{\sigma}^m}

\def\sm1{{\sigma}^{-1}}

\def\smtp1{{\sigma}^{-t+1}}

\def\o{\omega }
\def\S1{S^{-1}}

\def\Xpm1{X^{\pm 1}_1}

\def\sPM1{{\sigma }^{\pm 1}}
\def\sMP1{{\sigma }^{\mp 1 }}

\def\d{\delta}

\def\di{{\rm d.ind}}

\def\L{\Lambda}
\def\O{\Omega}

\def\G{\Gamma}

\def\CA{{\cal A}}

\def\Ytm1{Y^{t-1}}
\def\Yim1{Y^{i-1}}

\def\CK{{\cal K}}

\def\CM{{\cal M}}
\def\CN{{\cal N}}
\def\CS{{\cal S}}
\def\CF{{\cal F}}
\def\CG{{\cal G}}

\def\Aut{{\rm Aut}}

\def\Der{{\rm Der }}
\def\ad{{\rm ad }}

\def\ker{ {\rm ker } }

\def\SL2Z{ {\rm SL}_2({\bf Z}) }

\def\th{ \theta }

\def\Gp1{ G^{1 , 1 } }
\def\P11{ P^{-1 , 1 } }
\def\Pp1{ P^{1 , 1 } }

\def\th{\theta}

\def\nCLsr{{}^\nu\kern-2pt {\cal L}^{\sigma , \rho  }}
\def\nP{{}^\nu \kern-2pt P}
\def\nL{{}^\nu\kern-2pt L}
\def\nLL{{}^\nu\kern-2pt \Lambda}
\def\nPsr{{}^\nu\kern-2pt P^{\sigma , \rho  }}
\def\nLsr{{}^\nu\kern-2pt L^{\sigma , \rho  }}
\def\nuCL{{}^\nu\kern-2pt  {\cal L}}
\def\nCLsr{{}^\nu\kern-2pt {\cal L}^{\sigma , \rho  }}
\def\nCL1m{{}^\nu\kern-2pt {\cal L}^{-1 , 1  }}
\def\x1nu{x^\frac{1}{\nu}}
\def\xm1nu{x^{-\frac{1}{\nu}}}

\def\CN{{\cal N}}
\def\ra{\rightarrow }

\def\CB{{\cal B}}

\def\CI{{\cal I}}

\def\nAM0{{\nu }_{{\cal A},M_0}}
\def\nAN0{{\nu }_{{\cal A},N_0}}

\def\End{ {\rm End }}
\def\Der{ {\rm Der }}

\def\ad{ {\rm ad }}

\def\bx{\overline{x}}

\def\ga{\mathfrak{a}}

\def\gm{\mathfrak{m}}

\def\GL{{\rm GL}}
\def\SL{{\rm SL}}

\def\di!{\frac{\der^i}{i!}}
\def\dik!{\frac{\der^k_i}{k!}}

\def\N{\mathbb{N}}

\def\0{\overline{0}}
\def\1{\overline{1}}
\def\Lnz{\L_{n,\overline{0}}}
\def\Ln1{\L_{n,\overline{1}}}
\def\Lni{\L_{n,\overline{i}}}
\def\oa{\overline{a}}
\def\Lnod{\L_n^{od}}
\def\Lnev{\L_n^{ev}}
\def\Lnsod{\L_n'^{od}}
\def\Lnsev{\L_n'^{ev}}
\def\SDer{{\rm SDer}}

\def\Derod{\Der_K (\L_n)^{od}}
\def\Derev{\Der_K(\L_n)^{ev}}
\def\Derevs{\Der_K(\L_n)^{ev}_s}
\def\Derevs{\Der_K(\L_n)^{ev}_s}
\def\Derevn{\Der_K(\L_n)^{ev}_n}
\def\Derods{\Der_K(\L_n)^{od}_s}
\def\Derodn{\Der_K(\L_n)^{od}_n}
\def\SDerod{\SDer_K (\L_n)^{od}}
\def\SDerev{\SDer_K(\L_n)^{ev}}

\def\IDer{{\rm IDer}}
\def\ISDer{{\rm ISDer}}
\def\uod{u^{od}}
\def\uev{u^{ev}}
\def\dod{\d^{od}}
\def\dev{\d^{ev}}
\def\az{a_{\overline{0}}}
\def\a1{a_{\overline{1}}}
\def\sad{{\rm sad}}

\def\St{{\rm St}}
\def\S{\Sigma}

\def\CU{{\cal U}}
\def\Stab{{\rm Stab}}

\def\gl{{\rm gl}}
\def\DI{{\rm DI}}
\def\SDI{{\rm SDI}}
\def\hI{\widehat{I}}

\begin{document}

\author{V. V. \  Bavula 
}

\title{Derivations and skew derivations  of the Grassmann algebras}

\maketitle
\begin{abstract}
Surprisingly, skew derivations rather than ordinary derivations
are  more basic (important) object in study of the Grassmann
algebras. Let $\L_n = K\lfloor x_1, \ldots , x_n\rfloor$ be the
Grassmann algebra over a commutative ring $K$ with $\frac{1}{2}\in
K$, and $\d$ be a skew $K$-derivation of $\L_n$. It is proved that
$\d$ is a unique sum $\d = \d^{ev} +\d^{od}$ of an even and odd
skew derivation. Explicit formulae are given for $\d^{ev}$ and
$\d^{od}$ via the elements $\d (x_1), \ldots , \d (x_n)$. It is
proved that the set of all even skew derivations of $\L_n$
coincides with the set of all the inner skew derivations. Similar
results are proved for derivations of $\L_n$. In particular,
$\Der_K(\L_n)$ is a faithful but not simple $\Aut_K(\L_n)$-module
(where $K$ is reduced and $n\geq 2$). All differential and skew
differential ideals of $\L_n$ are found. It is proved that the set
of generic normal elements of $\L_n$  that are not units forms a
single $\Aut_K(\L_n)$-orbit (namely, $\Aut_K(\L_n)x_1$)   if $n$
is even and two orbits (namely, $\Aut_K(\L_n)x_1$ and
$\Aut_K(\L_n)(x_1+x_2\cdots x_n)$) if $n$ is odd.



$${\bf Contents}$$
\begin{enumerate}
\item Introduction. \item Derivations of the Grassmann rings.
\item Skew derivations of the Grassmann rings. \item Normal
elements of the Grassmann algebras.
\end{enumerate}
\end{abstract}


\section{Introduction}
Throughout, ring means an associative ring with $1$. Let $K$ be an
arbitrary ring (not necessarily commutative). The {\em Grassmann
algebra} (the {\em exterior algebra}) $\L_n = \L_n (K)= K\lfloor
x_1, \ldots , x_n\rfloor$ is generated freely over $K$ by elements
$x_1, \ldots , x_n$ that satisfy the defining relations:
$$ x_1^2=\cdots = x_n^2=0 \;\; {\rm and}\;\; x_ix_j=-x_jx_i\;\;
{\rm for \; all} \;\; i\neq j.$$ The Grassmann algebra $\L_n=
\oplus_{i\in \N}\L_{n,i}$ is an $\N$-graded algebra
($\L_{n,i}\L_{n,j} \subseteq \L_{n,i+j}$ for all $i,j\geq 0$)
where $\L_{n,i} := \oplus_{|\alpha | =i}Kx^\alpha$, $x^\alpha :=
x_1^{\alpha_1}\cdots x_n^{\alpha_n}$, and $|\alpha |:=
\alpha_1+\cdots +\alpha_n$.

{\bf Derivations of the Grassmann algebras}. Let $\Der_K(\L_n)$,
$\Der_K(\L_n)^{ev}$, $\Der_K(\L_n)^{od}$ and $\IDer_K(\L_n)$ be
the set of all, even, odd and inner derivations of $\L_n(K)$
respectively. Note that $\IDer_K(\L_n)= \{ \ad (a) \, | \, a\in
\L_n\}$ where $\ad (a) (x):= ax-xa$. Let $\Lnev$ and $\Lnod$ be
the set of even and odd elements of $\L_n$. Let
$\der_1:=\frac{\der}{\der x_1}, \ldots , \der_n:=\frac{\der}{\der
x_n}$ be partial skew $K$-derivations of $\L_n$ ($\der_i(x_j) =
\d_{ij}$, the Kronecker delta,  and $\der_i( a_ja_k) = \der_i
(a_j) a_k+(-1)^ja_j\der_i(a_k)$ for all $a_i\in \L_{n,i}$ and
$a_j\in \L_{n,j}$).

\begin{itemize}
\item (Theorem \ref{13Sep06}) {\em  Suppose that $K$ is a
commutative ring with $\frac{1}{2}\in K$. Then}
\begin{enumerate}
\item $\Der_K(\L_n) = \Derev \oplus \Derod$. \item $\Derev =
\oplus_{i=1}^n \Lnod \der_i$. \item $\Derod = \IDer_K(\L_n)$.
\item $\Der_K(\L_n)/ \IDer_K(\L_n)\simeq \Derev$.
\end{enumerate}
\end{itemize}
So, each derivation $\d \in \Der_K(\L_n)$ is a unique sum $\d =
\d^{ev} +\d^{od}$ of an even and odd derivation. When $K$ is a
field of characteristic $\neq 2$ this fact was proved by Djokovic,
 \cite{Djokovic78}. For an  even $n$,
let $ \Lnsod := \Lnod $. For an odd $n$, let $ \Lnsod$ be the
$K$-submodule of $\Lnod$ generated by all `monomials' $x^\alpha$
but $\th := x_1\cdots x_n$, i.e. $\Lnod = \Lnsod \oplus K\th$.
 The next result gives explicitly derivations $\d^{ev}$
and $\d^{od}$ via the elements $\d(x_1), \ldots , \d (x_n)$.

\begin{itemize}
\item (Corollary \ref{c13Sep06}) {\em Let $K$ be a commutative
ring with $\frac{1}{2}\in K$, $\d$ be a $K$-derivation of
$\L_n(K)$, and, for each $i=1, \ldots , n$, $\d (x_i) = \uev_i +
\uod_i$ for unique elements $\uev \in \Lnev $ and $\uod \in
\Lnod$. Then}
\begin{enumerate}
\item $\dev =\sum_{i=1}^n \uod_i\der_i$, {\em and} \item $\dod =
-\frac{1}{2}\ad (a)$ {\em where the {\em unique} element $a\in
\Lnsod$ is given by the formula}
$$ a= \sum_{i=1}^{n-1} x_1\cdots x_i\der_i\cdots \der_1\der_{i+1}
(\uev_{i+1})+\der_1 (\uev_1).$$
\end{enumerate}
\end{itemize}
The next results describes differential ideals of $\L_n$ (i.e.
which are stable under all derivations).
\begin{itemize}
\item (Proposition \ref{m19Sep06}) {\em Let $K$ be a commutative
ring with $\frac{1}{2}\in K$, $\CF_n(K):= \{ I: I_0\subseteq I_1
\subseteq \cdots \subseteq I_n\, | \, I_i\; {\rm are\;  ideals\;
of}\; K\}$ be the set of  $n$-flags of ideals of $K$, $\DI (\L_n)$
be the set of all differentiable ideals of $\L_n(K)$. Then the map
$$ \CF_n (K)\ra \DI (\L_n) , \; I\mapsto \hI :=
\oplus_{i=0}^n\oplus_{|\alpha | = i} I_ix^\alpha, $$ is a
bijection. In particular, $\gm^i$, $0\leq i\leq n+1$, are
differential ideals of $\L_n$; these are the only differential
ideals of $\L_n$ if $K$ is a field of characteristic} $\neq 2$.
\item (Theorem \ref{26Sep06}) {\em Let $K$ be a reduced
commutative ring with $\frac{1}{2}\in K$, $n\geq 1$. Then

1. $\Der_K(\L_n)$ is a faithful $\Aut_K(\L_n)$-module iff $n\geq
2$.

2. The $\Aut_K(\L_n)$-module $\Der_K(\L_n)$ is not simple. }
\end{itemize}

{\bf Skew derivations of the Grassmann algebras}. Let
$\SDer_K(\L_n)$, $\SDer_K(\L_n)^{ev}$, $\SDer_K(\L_n)^{od}$ and
$\ISDer_K(\L_n)$ be the set of all, even, odd and inner skew
derivations of $\L_n(K)$ respectively. $\ISDer_K(\L_n)=\{ \sad (a)
\, | \, a\in \L_n\}$ and $\sad (a) (a_i) := aa_i-(-1)^ia_ia$
($a_i\in \L_{n,i}$) is the inner skew derivation determined by the
element $a$. For an odd $n$, let $ \Lnsev := \Lnev $. For an even
$n$, let $ \Lnsev$ be the $K$-submodule of $\Lnev$ generated by
all `monomials' $x^\alpha$ but $\th := x_1\cdots x_n$, i.e. $\Lnev
= \Lnsev \oplus K\th$.

\begin{itemize}
\item (Theorem \ref{a13Sep06}) {\em Suppose that $K$ is a
commutative ring with $\frac{1}{2}\in K$. Then}
\begin{enumerate}
\item $\SDer_K(\L_n) = \SDerev \oplus \SDerod$. \item $\SDerod =
\oplus_{i=1}^n \Lnev\der_i$. \item $\SDerev = \ISDer_K(\L_n)$.
\item $\SDer_K(\L_n)/ \ISDer_K(\L_n)\simeq \SDerod$.
\end{enumerate}
\end{itemize}

So, any skew $K$-derivation $\d$ of $\L_n$ is a unique sum $\d =
\dev +\dod$ of an even and  odd
 skew derivation, and $\dev := \frac{1}{2}\sad (a)$ for a unique
element $a\in \Lnsev$. The next corollary describes explicitly the
skew derivations $\dev$ and $\dod$.

\begin{itemize}
\item (Corollary \ref{ca13Sep06})  {\em Let $K$ be a commutative
ring with $\frac{1}{2}\in K$, $\d$ be a skew $K$-derivation of
$\L_n(K)$, and, for each $i=1, \ldots , n$, $\d (x_i) = \uev_i +
\uod_i$ for unique  elements $\uev_i\in \Lnev $ and $\uod_i\in
\Lnod $. Then}
\begin{enumerate}
\item $\dod =\sum_{i=1}^n \uev_i\der_i$, {\em and} \item $\dev =
\frac{1}{2}\sad (a)$ {\em where the unique element $a\in \Lnsev$
is given by the formula}
$$ a= \sum_{i=1}^{n-1} x_1\cdots x_i\der_i\cdots \der_1\der_{i+1}
(\uod_{i+1})+\der_1 (\uod_1).$$
\end{enumerate}
\end{itemize}

{\bf The action of  $\Aut_K(\L_n)$ on the set of generic normal
non-units of $\L_n$}.

Let $\CN$ be the set of all the normal elements of the Grassmann
algebra $\L_n = \L_n(K)$,  $\CU$ be the set of all units of
$\L_n$, and $G:=\Aut_K(\L_n)$ be the group of $K$-automorphisms of
$\L_n$. Then $\CU \subseteq \CN$. The set $\CN$ is a disjoint
union of its $G$-invariant subsets,
$$\CN = \cup_{i=0}^n\CN_i, \;\; \CN_i:= \{ a\in \CN \, | \, a=
a_i+\cdots , 0\neq a_i\in \L_{n,i}\}.$$ Clearly, $\CN_0= \CU$. The
next result shows that `generic' normal non-unit elements of
$\L_n$ (i.e.  the set $\CN_1$) form a single $G$-orbit  if $n$ is
even, and two $G$-orbits if $n$ is odd.

\begin{itemize} \item  (Theorem \ref{5Nov06})  {\em  Let $K$ be a field of
characteristic $\neq 2$ and $\L_n = \L_n(K)$. Then}
\begin{enumerate} \item $\CN_1= Gx_1$ {\em if $n$ is
even.}\item $\CN_1= Gx_1\cup G(x_1+x_2\cdots x_n)$ {\em is the
disjoint union of two orbits if $n$ is odd.}
\end{enumerate}
\end{itemize}
The stabilizers of the elements $x_1$ and $x_1+x_2\cdots x_n$ are
found (Lemma \ref{c5Nov06} and Lemma \ref{s10Nov06}).


\section{Derivations of the Grassmann
rings}\label{DRGA}

In this section, the results on derivations from the Introduction
are proved. First, we recall some facts on Grassmann algebra (more
details the reader can find in \cite{BourbakiAlgCh1-3}).

{\bf The Grassmann algebra and its gradings}.  Let $K$ be an {\em
arbitrary} ring (not necessarily commutative). The {\em Grassmann
algebra} (the {\em exterior algebra}) $\L_n = \L_n (K)= K\lfloor
x_1, \ldots , x_n\rfloor$ is generated freely over $K$ by elements
$x_1, \ldots , x_n$ that satisfy the defining relations:
$$ x_1^2=\cdots = x_n^2=0 \;\; {\rm and}\;\; x_ix_j=-x_jx_i\;\;
{\rm for \; all} \;\; i\neq j.$$ Let $\CB_n$ be the set of all
subsets of the set of indices $\{ 1, \ldots , n\}$. We may
identify the set $\CB_n$ with the direct product $\{ 0,1\}^n$ of
$n$ copies of the two-element set $\{ 0, 1\}$ by the rule $\{ i_1,
\ldots , i_k\} \mapsto (0, \ldots, 1, \ldots , 1, \ldots , 0)$
where $1$'s are on $i_1, \ldots , i_k$ places and $0$'s elsewhere.
So, the set $\{ 0, 1\}^n$ is the set of all the characteristic
functions on the set $\{ 1, \ldots , n\}$.
$$ \L_n = \bigoplus_{\alpha \in \CB_n} Kx^\alpha = \bigoplus_{\alpha \in \CB_n} x^\alpha K,
 \;\; x^\alpha :=x_1^{\alpha_1} \cdots  x_n^{\alpha_n}, $$
where $\alpha = (\alpha_1, \ldots , \alpha_n)\in \{ 0, 1\}^n =
\CB_n$.  Note that the order in the product $x^\alpha$ is fixed.
So, $\L_n$ is a free left and right $K$-module of rank $2^n$. Note
that $(x_i):= x_i\L_n = \L_nx_i$ is an ideal of $\L_n$. Each
element $a\in \L_n$ is a unique sum $a= \sum a_\alpha x^\alpha$,
$a_\alpha \in K$. One can view each element $a$ of $\L_n$ as a
`function' $a= a(x_1, \ldots , x_n)$ in the non-commutative
variables $x_i$. The $K$-algebra epimorphism
\begin{eqnarray*}
 \L_n &\ra & \L_n / (x_{i_1}, \ldots , x_{i_k})\simeq K\lfloor x_1, \ldots
, \widehat{x_{i_1}}, \ldots , \widehat{x_{i_k}}, \ldots ,
x_n\rfloor, \\
a&\mapsto & a|_{x_{i_1}=0, \ldots , x_{i_k}=0} :=a+ (x_{i_1},
\ldots , x_{i_k}),
\end{eqnarray*}
 may be seen as the operation of
taking value of the function $a(x_1, \ldots , x_n)$ at the point
$x_{i_1} = \cdots = x_{i_k}=0$ where here and later the hat over a
symbol means that it is missed.

For each $\alpha \in \CB_n$, let $ | \alpha | := \alpha_1 +\cdots
+ \alpha_n$. The ring $\L_n=\oplus_{i=0}^n \L_{n,i}$ is a
$\Z$-{\em graded} ring ($\L_{n,i}\L_{n,j}\subseteq \L_{n,i+j}$ for
all $i,j$) where $\L_{n,i}:= \oplus_{|\alpha | =i}Kx^\alpha$. The
ideal $\gm := \oplus_{i\geq 1} \L_{n,i}$ of $\L_n$  is called the
{\em augmentation} ideal. Clearly, $K\simeq \L_n/ \gm$, $\gm^n=
Kx_1\cdots x_n$ and $\gm^{n+1}=0$. We say that an element $\alpha$
of $\CB_n$ is {\em even} (resp. {\em odd}) if the set $\alpha$
contains even (resp. odd) number of elements. By definition, the
empty set is even. Let $\Z_2:= \Z / 2\Z = \{\0, \1 \}$. The ring
$\L_n = \Lnz \oplus \Ln1$ is a $\Z_2$-{\em graded} ring where
$\Lnz := \Lnev :=\oplus_{\alpha \; {\rm is \; even}}Kx^\alpha$ is
the subring of even elements of $\L_n$ and $\Ln1 := \Lnod
:=\oplus_{\alpha \; {\rm is \; odd}}Kx^\alpha$ is the
$\Lnev$-module of odd elements of $\L_n$. The ring $\L_n$ has the
$\gm$-{\em adic} filtration $\{ \gm^i\}_{i\geq 0}$. The even
subring $\Lnev$ has the induced $\gm$-adic filtration $\{ \L_{n,
\geq i}^{ev} := \Lnev \cap \gm^i \}$. The $\Lnev$-module $\Lnod$
has the induced $\gm$-adic filtration $\{ \L_{n, \geq i}^{od} :=
\Lnod \cap \gm^i \}$.

 The
$K$-linear map $ a\mapsto \oa$ from $\L_n$ to itself which is
given by the rule
$$\oa :=
\begin{cases}
a,& \text{if $a\in \Lnz$},\\
-a,& \text{if $a\in \Ln1$},
\end{cases}$$
is a ring automorphism such that  $\overline{\oa}= a$ for all
$a\in \L_n$. For all $a\in \L_n$ and $i=1, \ldots , n$,
\begin{equation}\label{xiaib}
x_ia= \oa x_i \;\; {\rm and}\;\; ax_i= x_i\oa .
\end{equation}
So, each element $x_i$ of $\L_n$ is a {\em normal} element, i.e.
the two-sided ideal $(x_i)$ generated by the element $x_i$
coincides with both left and right ideals generated by $x_i$:
$(x_i) = \L_n x_i= x_i\L_n$.

For an arbitrary $\Z$-graded ring $A= \oplus_{i\in \Z} A_i$,  an
additive map $\d :A\ra A$ is called a {\em left skew derivation}
if 
\begin{equation}\label{dsd}
\d (a_ia_j) = \d (a_i) a_j+(-1)^ia_i\d (a_j)\;\; {\rm for \;
all}\;\; a_i\in A_i, \; a_j\in A_j.
\end{equation}
In this paper, a skew derivation means a {\em left} skew
derivation. Clearly, $1\in \ker (\d )$ ($ \d (1) = \d (1\cdot 1) =
2\d (1)$ and so $\d (1)=0$). The restriction of the left skew
derivation $\d$ to the even subring $A^{ev}:= \oplus_{i\in 2\Z}
A_i$ of $A$ is an ordinary derivation. Recall that an additive
subgroup $B$ of $A$ is called a homogeneous subgroup if $B=
\oplus_{i\in \Z } B\cap A_i$.

{\it Definition}. For the  ring $\L_n (K)$, consider the set of
{\em left} skew $K$-derivations:
$$ \der_1:= \frac{\der}{\der x_1}, \ldots , \der_n:=
\frac{\der}{\der x_n}$$ given by the rule $\der_i (x_j)= \d_{ij}$,
the Kronecker delta. Informally, these skew $K$-derivations will
be called (left) {\em partial skew derivatives}.

{\it Example}. $\der_i (x_1 \cdots x_i \cdots x_k)= (-1)^{i-1}
x_1\cdots x_{i-1} x_{i+1} \cdots x_k$.

If  the ring $K$ is commutative, $2\in K$ is regular (i.e. $2\l
=0$ in $K$ implies $\l =0$), and $n\geq 2$,  then the centre of
$\L_n(K)$ is equal to
$$ Z(\L_n)=\begin{cases}
\Lnz , & \text{if $n$ is even},\\
\Lnz\oplus Kx_1\cdots x_n, & \text{if $n$ is odd }.
\end{cases}$$

Let $K$ be a commutative ring. For an  even $n$, let $ \Lnsod :=
\Lnod $. For an odd $n$, let $ \Lnsod$ be the $K$-submodule of
$\Lnod$ generated by all `monomials' $x^\alpha$ but $\th :=
x_1\cdots x_n$, i.e. $\Lnod = \Lnsod \oplus K\th$.
 Similarly, for an odd $n$, let $
\Lnsev := \Lnev $. For an even $n$, let $ \Lnsev$ be the
$K$-submodule of $\Lnev$ generated by all `monomials' $x^\alpha$
but $\th := x_1\cdots x_n$, i.e. $\Lnev = \Lnsev \oplus K\th$.
 For any $n$,
\begin{equation}\label{Lndd}
\Lnod = \Lnsod \oplus \Lnod\cap Z(\L_n ).
\end{equation}
So, one can naturally identify $\Lnod / \Lnod \cap Z(\L_n )$ with
$\Lnsod$.

Consider the sets of {\em even} and {\em odd} $K$-derivations of
$\L_n(K)$:
\begin{eqnarray*}
\Derev & :=& \{ \d \in \Der_K(\L_n )\, | \, \d (\Lni ) \subseteq \Lni , \; \overline{i} \in \Z_2 \},\\
 \Derod & :=& \{ \d \in \Der_K(\L_n )\, | \, \d (\Lni ) \subseteq \L_{n, \overline{i}+\overline{1}} , \; \overline{i} \in \Z_2
 \}.
\end{eqnarray*}
So, even derivations are precisely the derivations that respect
$\Z_2$-grading of $\L_n$, and the odd derivations are precisely
the derivations that reverse it. The set of odd and even
derivations are left $\Lnz$-modules. For each element $a\in \L_n$,
one can attach the $K$-derivation of $\L_n$ $\ad (a) : b\mapsto
[a,b]:= ab-ba$, so-called, the {\em inner} derivation determined
by $a$. The set of all inner derivations is denoted by
$\IDer_K(\L_n)$, and the map
$$\L_n/ Z(\L_n)\ra \IDer_K(\L_n), \;\; a+Z(\L_n) \mapsto \ad (a),
$$
 is an isomorphism of left  $Z(\L_n)$-modules. The next theorem describes
 explicitly the sets of all/inner/even and odd derivations.
\begin{theorem}\label{13Sep06}
Suppose that $K$ is a commutative ring with $\frac{1}{2}\in K$.
Then
\begin{enumerate}
\item $\Der_K(\L_n) = \Derev \oplus \Derod$. \item $\Derev =
\oplus_{i=1}^n \Lnod \der_i$. \item $\Derod = \IDer_K(\L_n)$ and
the map
$$ \ad : \Lnsod = \Lnod / \Lnod \cap Z(\L_n )\ra \IDer_K(\L_n),
\;\; a\mapsto \ad (a), $$ is the $Z(\L_n)$-module isomorphism.
\item $\Der_K(\L_n)/ \IDer_K(\L_n)\simeq \Derev$.
\end{enumerate}
\end{theorem}

{\it Proof}. Since $\Derev \cap \Derod = 0$, one has the inclusion
\begin{equation}\label{1ri}
\Der_K(\L_n)\supseteq \Derev \oplus \Derod .
\end{equation}
Clearly, 
\begin{equation}\label{2ri}
\Derev \supseteq \sum_{i=1}^n \Lnod \der_i =\bigoplus_{i=1}^n
\Lnod \der_i ,
\end{equation}
$\IDer_K(\L_n) \simeq \L_n / Z(\L_n) = (\Lnev \oplus \Lnod ) /
Z(\L_n)\simeq \Lnod / \Lnod \cap Z(\L_n)\simeq \Lnsod$ since
$\Lnev \subseteq Z(\L_n)$. For each $a\in \Lnod$, $\ad (a) \in
\Derod$, hence 
\begin{equation}\label{3ri}
\Derod \supseteq \IDer_K(\L_n).
\end{equation}
Note that statement 4 follows from statements 1 and 3. Now, it is
obvious that in order to finish the proof of the theorem it
suffices to show that

{\it Claim}. $\Der_K(\L_n) \subseteq \sum_{i=1}^n \Lnod \der_i
+\IDer_K(\L_n)$.

Indeed, suppose that the inclusion of the claim holds then, by
(\ref{2ri}) and (\ref{3ri}),
$$ \Der_K(\L_n) \subseteq \sum_{i=1}^n \Lnod \der_i
+\IDer_K(\L_n)\subseteq \Derev \oplus \Derod,$$ hence statement 1
is true by (\ref{1ri}). Statement 1 together with inclusions
(\ref{2ri}) and (\ref{3ri}) implies statements 2 and 3.

{\it Proof of the Claim}. Let $\d$ be a $K$-derivation of $\L_n$.
We have to represent the derivation $\d $ as a sum
$$ \d = \sum_{i=1}^n a_i\der_i + \ad (a), \;\; a_i\in \Lnod, \;\;
a\in \L_n.$$ The proof of the claim is constructive.
 According to the decomposition $\L_n= \Lnev \oplus \Lnod$ each
element $u$ of $\L_n$ is a unique sum 
\begin{equation}\label{u=ueo}
u= \uev + \uod
\end{equation}
 of its even and odd components ($\uev \in \Lnev$
and $\uod \in \Lnod$). For each $i$, let $u_i:= \d (x_i) =
\uev_i+\uod_i$,
$$\der := \sum_{i=1}^n \uod_i\der_i \;\; {\rm and}\;\; \d':= \d -
\der.$$ Note that $\der \in \sum_{i=1}^n \Lnod \der_i$, hence
changing $\d$ for $\d'$, if necessary, one may assume that  all
the  elements  $u_i$ are even. So,  it suffices to show that $\d =
\ad (a)$ for some $a$. We produce such an $a$ in several steps.

{\it Step 1}. Let us prove that, {\em for each $i=1, \ldots , n$,
$u_i= v_ix_i$ for some element } $v_i\in K\lfloor x_1, \ldots ,
\widehat{x_i}, \ldots , x_n\rfloor^{od}$. Note that $0= \d (0)= \d
(x_i^2) = u_ix_i + x_iu_i= 2u_ix_i$, and so $u_ix_i=0$ (since
$\frac{1}{2}\in K$). This means that $u_i= v_ix_i$ for some
element $v_i\in K\lfloor x_1, \ldots , \widehat{x_i}, \ldots ,
x_n\rfloor^{od}$ since $u_i$ is even. For $n=1$, it gives $u_1=0$
 since $K^{od}=0$, and we are done. So, let $n\geq 2$.

{\it Step 2}. We claim that, {\em for each  pair} $i\neq j$,
\begin{equation}\label{ast2}
v_i|_{x_j=0}=v_j|_{x_i=0}.
\end{equation}
Evaluating the derivation $\d$ at the element $0= x_ix_j+x_jx_i$
and taking into account that all the elements $u_i$ are even
(hence central) we obtain
$$ 0= 2( u_ix_j+u_jx_i) = 2(v_ix_ix_j+v_jx_jx_i) = 2( v_i-v_j)
x_ix_j.$$ This means that $v_i-v_j\in (x_i, x_j)$ since
$\frac{1}{2}\in K$, or, equivalently, $v_i|_{x_i=0,
x_j=0}=v_j|_{x_i=0, x_j=0}$. By Step 1, this equality can be
written as (\ref{ast2}).

{\it Step 3}. Note that $\d (x_1)=v_1x_1= \ad (\frac{1}{2}v_1)
(x_1)$ since  $v_1$ is odd, i.e. $(\d -\ad
(\frac{1}{2}v_1))(x_1)=0$. So, changing $\d$ for $\d -\ad
(\frac{1}{2}v_1)$ one can assume that $\d (x_1)=0$, i.e. $v_1=0$.
Then, by (\ref{ast2}), $v_i|_{x_1=0}=0$ for all $i=2, \ldots , n$,
and so $v_i\in (x_1 x_i)$ for all $i=2, \ldots , n$. Summarizing,
we can say  that by adding to $\d$ a well chosen  inner derivation
one can assume that $\d (x_1)=0$ and $\d (x_i) \in (x_1 x_i)$ for
all $i\geq 2$. This statement serves as the base of the induction
in the proof  of the next statement. For each $k$ such that $1\leq
k\leq n$, by adding to $\d$ a  certain  inner derivation we can
assume that 
\begin{equation}\label{kind}
\d (x_1)=\cdots = \d (x_k)=0, \;\; \d (x_i) \in (x_1 \cdots
x_kx_i), \;\; k<i\leq n .
\end{equation}
So, assuming that (\ref{kind}) holds for $k$ we must prove  the
same statement but for $k+1$. Note that $v_{k+1} x_{k+1} = \d (
x_{k+1})\in (x_1\cdots x_k)$, hence $v_{k+1} = x_1\cdots x_kv$ for
some $v\in K\lfloor x_{k+2} , \ldots , x_n\rfloor$. Consider the
derivation $\d':= \d -\ad (\frac{1}{2}v_{k+1})$. For each $i=1,
\ldots , k$, $\d'(x_i)= \d (x_i) =0$ as $v_{k+1} \in ( x_1\cdots
x_k)$; and $\d'( x_{k+1}) = v_{k+1} x_{k+1}- [\frac{1}{2}v_{k+1},
x_{k+1}]= v_{k+1} x_{k+1} - v_{k+1} x_{k+1}=0$. These prove the
first part of (\ref{kind}) for $k+1$, namely, that
$$ \d' (x_1)=\cdots = \d' (x_{k+1})=0.$$
So, changing $\d$ for $\d'$ one can assume that
$$ \d (x_1)=\cdots = \d (x_{k+1})=0.$$
These conditions imply that $v_1=\cdots = v_{k+1}=0$. If $n=k+1$,
we are done. So, let $k+1<n$. Then, by (\ref{ast2}), for each
$i>k+1$, $v_i \in \cap_{j=1}^{k+1} (x_j)= (x_1\cdots x_{k+1})$,
hence $\d (x_i) = v_ix_i \in ( x_1\cdots x_{k+1} x_i)$. By
induction, (\ref{kind}) is true for all $k$. In particular, for
$k=n$ one has $\d =0$. This means that $\d$ is an inner
derivation, as required. $\Box $

 The ring $K[x]/(x^2)$ of dual
numbers is the Grassmann ring $\L_1$.

\begin{corollary}\label{cc13Sep06}
Suppose that $K$ is a commutative ring with $\frac{1}{2}\in K$.
Then $\Der_K(K[x]/(x^2))=\Der_K(K[x]/(x^2))^{ev} = Kx\frac{d}{dx}$
and $\Der_K(K[x]/(x^2))^{od}=0$  where $\frac{d}{dx}$ is the skew
$K$-derivation of $K[x]/(x^2)$.
\end{corollary}

A Lie algebra $(\CG , [\cdot , \cdot ])$ over $K$ is  positively
graded  if $\CG = \oplus_{i\geq 0} \CG_i$ is a direct sum of
$K$-submodules such that $[\CG_i, \CG_j]\subseteq \CG_{i+j}$ for
all $i,j\geq 0$.

$(\Der_K(\L_n), [\cdot , \cdot ])$ is a Lie algebra over $K$ where
$[\d , \der ] := \d \der - \der \d$. By Theorem \ref{13Sep06}, the
Lie algebra $\Der_K(\L_n)=\oplus_{i\geq 0}D_i$ is a positively
graded Lie algebra where $$D_i:= \{ \d \in \Der_K(\L_n)\, | \, \d
(\L_{n,j})\subseteq \L_{n,j+i}, \; j\geq 0\}.$$ Clearly, $D_i=0$,
$i\geq n$. For each even natural number $i$ such that $0\leq i\leq
n-1$, 
\begin{equation}\label{Diev}
D_i=\oplus_{j=1}^n\L_{n,i+1}\der_j.
\end{equation}
For each odd natural number $i$ such that $1\leq i\leq n-1$,
\begin{equation}\label{Diod}
D_i=\{ \ad (a) \, | \, a\in \L_{n,i}\} \simeq \L_{n, i},\;\;  \ad
(a) \mapsto a.
\end{equation}
The zero component $D_0= \oplus_{i,j=0}^nKx_i\der_j$ of
$\Der_K(\L_n)$ is a Lie subalgebra of $\Der_K(\L_n)$ which is
canonically  isomorphic to the Lie algebra $\gl_n(K):=
\oplus_{i,j=1}^nKE_{ij}$ via $D_0\ra \gl_n(K)$, $x_i\der_j\mapsto
E_{ij}$, where $E_{ij}$ are the matrix units.  By the very
definition, $D_+:= \oplus_{i\geq 1}D_i$ is a nilpotent ideal of
the Lie algebra $\Der_k(\L_n)$ such that $\Der_K(\L_n)= D_0\oplus
D_+\simeq \gl_n(K)\oplus D_+$ and $\Der_K(\L_n)/D_+\simeq
\gl_n(K)$. So, if $K$ is a field of characteristic zero then $D_+$
is the radical of the Lie algebra $\Der_K(\L_n)$.

The Lie algebra $\Der_K(\L_n)= \Derev \oplus \Derod$ is a
$\Z_2$-graded Lie algebra.

By Theorem \ref{13Sep06}, any $K$-derivation $\d$ of $\L_n$ is a
unique sum $\d = \dev +\dod$ of even and odd derivations, and
$\dod := -\frac{1}{2}\ad (a)$ for a {\em unique} element $a\in
\Lnsod$. In order to find the element $a$ (Corollary
\ref{c13Sep06}), we need two theorems which are interesting on
their own right. Theorem \ref{14Sep06} gives a unique (sort of
`triangular') canonical presentation of any element of $\L_n$.
This presentation is important in dealing with derivations and
skew derivations. The element $a$ in $\dod = -\frac{1}{2}\ad (a)$
 is given in this  form (Corollary \ref{c13Sep06}). In order to find the element
$a$ we need to find solutions to the system of equations (Theorem
\ref{s14Sep06}). This system is a kind of Poincar\'{e} Lemma for
the  (noncommutative) Grassmann algebra $\L_n$.

\begin{theorem}\label{14Sep06}
\cite{jacgras} Let $K$ be an arbitrary (not necessarily
commutative) ring. Then
\begin{enumerate}
\item the Grassmann ring $\L_n(K)$ is a direct sum of right
$K$-modules
\begin{eqnarray*}
\L_n(K)&=& x_1\cdots x_nK \oplus x_1\cdots x_{n-1}K \oplus
x_1\cdots x_{n-2}K\lfloor x_n\rfloor \oplus\cdots \\
&\cdots &\oplus x_1\cdots x_iK\lfloor x_{i+2}\ldots ,  x_n\rfloor
\oplus\cdots \oplus  x_1 K\lfloor x_3\ldots ,  x_n\rfloor\oplus
K\lfloor x_2\ldots , x_n\rfloor .
\end{eqnarray*}
\item So, each element $a\in \L_n(K)$ is a unique sum
$$ a= x_1\cdots x_na_n+ x_1\cdots x_{n-1}b_n+\sum_{i=1}^{n-2}
x_1\cdots x_ib_{i+1} + b_1$$ where $a_n, b_n\in K$, $b_i\in
 K\lfloor x_{i+1}\ldots ,  x_n\rfloor$, $1\leq i\leq n-1$.
 Moreover,
 \begin{eqnarray*}
 a_n&=& \der_n\der_{n-1}\cdots \der_1(a), \\
b_{i+1}&=&\der_i\der_{i-1}\cdots \der_1(1-x_{i+1}\der_{i+1})(a), \; 1\leq i\leq n-1, \\
 b_1&=&(1-x_1\der_1)(a).
\end{eqnarray*}
So, $$ a= x_1\cdots x_n\der_n\der_{n-1}\cdots
\der_1(a)+\sum_{i=1}^{n-1} x_1\cdots x_i\der_i\cdots
\der_1(1-x_{i+1}\der_{i+1})(a)+(1-x_1\der_1)(a).$$
\end{enumerate}
\end{theorem}

By Theorem \ref{14Sep06}, the identity map ${\rm id}_{\L_n}
:\L_n\ra \L_n$ is equal to 
\begin{equation}\label{idLn}
{\rm id}_{\L_n}=x_1\cdots x_n\der_n\der_{n-1}\cdots \der_1+
\sum_{i=1}^{n-1} x_1\cdots x_i\der_i\cdots
\der_1(1-x_{i+1}\der_{i+1})+(1-x_1\der_1).
\end{equation}

\begin{theorem}\label{s14Sep06}
\cite{jacgras} Let $K$ be an arbitrary ring, $u_1, \ldots , u_n\in
\L_n(K)$, and $ a\in \L_n(K)$ be an unknown. Then the system of
equations
$$\begin{cases}
x_1a=u_1 \\
x_2a=u_2 \\
\;\;\;\; \;\;\;\vdots \\
x_na=u_n
\end{cases}
$$
has a solution in $\L_n$ iff the following two conditions hold
\begin{enumerate}
\item $u_1\in (x_1), \ldots , u_n\in (x_n)$, and \item
$x_iu_j=-x_ju_i$ for all $i\neq j$.
\end{enumerate}
In this case, 
\begin{equation}\label{alsol}
a= x_1\cdots x_na_n+\sum_{i=1}^{n-1} x_1\cdots x_i\der_i\cdots
\der_1\der_{i+1}(u_{i+1})+\der_1(u_1), \;\; a_n \in K,
\end{equation}
are all the solutions.
\end{theorem}

The next corollary describes explicitly $\dev$ and $\dod$ in $\d =
\dev + \dod$.

\begin{corollary}\label{c13Sep06}
Let $K$ be a commutative ring with $\frac{1}{2}\in K$, $\d$ be a
$K$-derivation of $\L_n(K)$, and, for each $i=1, \ldots , n$,  $\d
(x_i) = \uev_i + \uod_i$ for unique  elements $\uev \in \Lnev $
and $\uod \in \Lnod$. Then
\begin{enumerate}
\item $\dev =\sum_{i=1}^n \uod_i\der_i$, and \item $\dod =
-\frac{1}{2}\ad (a)$ where the {\em unique} element $a\in \Lnsod$
is given by the formula
$$ a= \sum_{i=1}^{n-1} x_1\cdots x_i\der_i\cdots \der_1\der_{i+1}
(\uev_{i+1})+\der_1 (\uev_1).$$
\end{enumerate}
\end{corollary}

{\it Proof}. 1. This statement  has  been proved already in the
proof of Theorem \ref{13Sep06}.

2. For each $i=1, \ldots , n$, on the one hand $\dod (x_i) = (\d -
\dev )(x_i) = \uev_i$; on the other, $ \dev (x_i) =
-\frac{1}{2}(ax_i-x_ia) = \frac{1}{2}2x_ia= x_ia$. So, the element
$a$ is a solution to the system of  equations
$$\begin{cases}
x_1a=\uev_1 \\
x_2a=\uev_2 \\
\;\;\;\; \;\;\;\vdots \\
x_na=\uev_n.
\end{cases}
$$
By Theorem \ref{s14Sep06} and the fact that $a\in \Lnsod$ (i.e.
$a_n=0$), we have
$$ a= \sum_{i=1}^{n-1} x_1\cdots x_i\der_i\cdots \der_1\der_{i+1}
(\uev_{i+1})+\der_1 (\uev_1).\;\;\; \Box$$

Let $K$ be a field. Let $V$ be a finite dimensional vector space
over $K$ and $a\in \End_K(V)$, a $K$-linear map on $V$. The vector
space $V$ is the $K[t]$-module where $t\cdot v := av$; $V$ is the
$K[a]$-module for short. The linear map $a$ is called {\em
semi-simple} (resp. {\em nilpotent}) if the $K[a]$-module $V$ is
semi-simple (resp. $a^k=0$ for some $k\geq 1$). It is well-known
that $a$ is a unique sum $a= a_s+a_n$ where $a_s$ is a semi-simple
map, $a_n$ is a nilpotent map, and $a_s, a_n\in K[a]:= \sum_{i\geq
0} Ka^i$ (in particular, the  maps $a$, $a_s$, and $a_n$ commute).
If $V$ is a finite dimensional $K$-algebra and $a$ is a
$K$-derivation of the algebra $V$ then the maps $a_s$ and $a_n$
are also $K$-derivations.

The subsets of $\Der_K(V)$ of all semi-simple derivations
$\Der_K(V)_s$ and all nilpotent derivations  $\Der_K(V)_n$ do not
meet, i.e. $\Der_K(V)_s\cap \Der_K(V)_n=0$. In general, the sets
of semi-simple and nilpotent derivations are {\em not} vector
spaces, though they are closed under scalar multiplication.

Let $K$ be a commutative ring with $\frac{1}{2}\in K$. The
following $K$-derivations of $\L_n= \L_n(K)$
$$ h_1:= x_1\der_1, \ldots , h_n:= x_n\der_n,$$
commute, $h_1^2=h_1, \ldots , h_n^2=h_n$, and $h_i(x^\alpha ) =
\alpha_ix^\alpha$ for all $i$ and $\alpha$. Since $h_i(x_j)=
\d_{ij}x_j$,  the maps $h_1, \ldots , h_n$ are linearly
independent. So, $h_1, \ldots , h_n$ are {\em commuting,
semi-simple, $K$-linearly independent, idempotent $K$-derivations}
of the algebra $\L_n$. For each $i=1, \ldots , n$, $\L_n=
K_i\oplus x_iK_i$ where $K_i:= \ker (\der_i) = K\lfloor x_1,
\ldots , \widehat{x_i}, \ldots , x_n\rfloor$, and $h_i: \L_n\ra
\L_n$ is the projection onto $x_iK_i$.

Let $H$ be the subalgebra of the endomorphism algebra
$\End_K(\L_n)$ generated by the elements $h_1, \ldots , h_n$. As
an abstract algebra $H\simeq K[H_1, \ldots , H_n]/ (H_1^2, \ldots
, H_n^2)$. The algebra $\L_n=\oplus_{\alpha \in \CB_n} Kx^\alpha$
is a semi-simple $H$-module where each isotypic component is
simple: $Kx^\alpha$ is the simple $H$-module, and $Kx^\alpha\simeq
Kx^\beta$ as $H$-modules  iff $\alpha = \beta$. Let
\begin{eqnarray*}
\Derevs:= \Der_K(\L_n)^{ev}\cap \Der(\L_n)_s &
\Derevn :=\Der_K(\L_n)^{ev}\cap \Der(\L_n)_n,\\
\Derods:= \Der_K(\L_n)^{od}\cap \Der(\L_n)_s, & \Derodn
:=\Der_K(\L_n)^{od}\cap \Der(\L_n)_n.
\end{eqnarray*}

{\it Definition}.  An ideal $\ga$ of $\L_n$ is called a {\em
differential} ideal (or a $\Der_K(\L_n)$-{\em invariant} ideal) if
$\d (\ga ) \subseteq \ga$ for all $\d \in \Der_K(\L_n)$.

The next proposition describes all the differential ideals of
$\L_n(K)$.

\begin{proposition}\label{m19Sep06}
Let $K$ be a commutative ring with $\frac{1}{2}\in K$, $\CF_n(K):=
\{ I: I_0\subseteq I_1 \subseteq \cdots \subseteq I_n\, | \, I_i\;
{\rm are\;  ideals\; of}\; K\}$ be the set of  $n$-flags of ideals
of $K$, $\DI (\L_n)$ be the set of all differentiable ideals of
$\L_n(K)$. Then the map
$$ \CF_n (K)\ra \DI (\L_n) , \; I\mapsto \hI :=
\oplus_{i=0}^n\oplus_{|\alpha | = i} I_ix^\alpha, $$ is a
bijection. In particular, $\gm^i$, $0\leq i\leq n+1$, are
differential ideals of $\L_n$; these are the only differential
ideals of $\L_n$ if $K$ is a field of characteristic $\neq 2$.
\end{proposition}

{\it Proof}. Recall that $\Der_K(\L_n) =\oplus_{i\geq 0} D_i$. By
(\ref{Diev}) and (\ref{Diod}), the map $I\mapsto \hI$ is
well-defined and injective, by the very definition. It remains to
show that each differential idea, say $\ga$, of $\L_n(K)$ is equal
to $\hI$ for some $I\in \CF_n(K)$. Since $\sum_{i=1}^n
Kh_i\subseteq \Der_K(\L_n)$ and $ \L_n =\oplus_{\alpha \in \CB_n}
Kx^\alpha$ is the direct sum of non-isomorphic simple $H$-modules,
we have
 $$\ga = \oplus_{\alpha \in
\CB_n}(\ga \cap Kx^\alpha )= \oplus_{\alpha \in \CB_n} \ga_\alpha
x^\alpha$$ where $\ga_\alpha$ is an ideal of $K$ such that
$\ga_\alpha x^\alpha = \ga \cap Kx^\alpha$. For each $i$ such that
$0\leq i \leq n$, $\oplus_{|\alpha |=i}\ga_\alpha x^\alpha$ is a
$D_0$-module where $D_0=\oplus_{i,j=1}^n Kx_i\der_j$, hence all
the ideals $\ga_\alpha$ coincide where $|\alpha | = i$. Let $I_i$
be their common value. Since $\ga$ is an ideal of $\L_n$, $\{
I:I_0\subseteq \cdots \subseteq I_n\} \in \CF_n(K)$, and so $\ga =
\hI$, as required. $\Box$

{\it Definition}. A ring $R$ is called a {\em differentiably
simple} ring if it is a simple left $\Der (R)$-module.

So, the algebra $\L_n$ is {\em not} differentiably simple if
$n\geq 1$ where $K$ is a commutative ring with $\frac{1}{2}\in K$.

Let $K$ be a {\em reduced} commutative ring with $\frac{1}{2}\in
K$. Let $\CS$ be the set of all $n$-tuples $(s_1, \ldots , s_n)$
where $s_1, \ldots , s_n$ are commuting, idempotent
$K$-derivations (i.e. $s_i^2=s_i$) of $\L_n$ such that the
following conditions hold: $\cap_{i=1}^n \ker (s_i) = K$;  all the
$K$-modules $\CK_i:= \gm \cap \ker (s_i-1) \cap \cap_{j\neq i}
\ker (s_j)$ are  free of rank $1$ over $K$, i.e. $\CK_i =
Kx_i'\simeq K$ for some element $x_i'\in \gm$ $K$; $\L_n = \ker
(s_1)+\CK_1\ker (s_1)$; and, for each $i=1, \ldots , n-1$,
\begin{equation}\label{K1i}
K_{1, \ldots , i}= K_{1, \ldots , i+1}+\CK_{i+1} K_{1, \ldots ,
i+1}
\end{equation}
where $K_{1, \ldots , i}:= \cap_{j=1}^i \ker(s_j)$.

Clearly, $(h_1, \ldots, h_n)\in \CS$ (see (\ref{CKxi}) below). For
each $(s_1, \ldots, s_n)\in \CS$, the maps $s_1, \ldots , s_n$ are
$K$-linearly independent ($\sum \mu_i s_i=0$ $\Rightarrow$
$0=(\sum \mu_is_i)(\CK_i)= \mu_iKx_i'$ $\Rightarrow$ $\mu_i=0$).

Let $G:= \Aut_K(\L_n)$ be the group of $K$-algebra automorphisms
of the Grassmann algebra $\L_n$.  For $\s \in G$, let $x_i':= \s
(x_i)$. Then $x_i'^2=\s(x_i^2)= \s (0)=0$. If $ \l_i\equiv
x_i'\mod \gm$ for some $\l_i\in K$ then $\l_i^2=0$, hence $\l_i=0$
since $K$ is reduced. Therefore, $\s (\gm ) = \gm$, and so
\begin{equation}\label{smi=mi}
\s (\gm^i)=\gm^i \;\; {\rm for \; all}\;\;  i\geq 1.
\end{equation}
  By (\ref{smi=mi}), the group $G$
acts on the set $\CS$ by conjugation (i.e. by changing
generators): $\s \cdot (s_1, \ldots, s_n):= (\s s_1\s^{-1},
\ldots, \s s_n\s^{-1})$ We prove shortly that the group $G$ act
transitively on the set $\CS$ (Corollary \ref{a19Sep06}.(2)) and
the stabilizer $\St (h_1, \ldots , h_n)$ of the element $(h_1,
\ldots , h_n)$ is equal to the `$n$-{\em dimensional algebraic
torus}' (where $K^*$ is the group of units of $K$)
$$ \mathbb{T}^n:= \{ \s_\l \, | \, \l \in K^{*n}, \s_\l ( x_i) =
\l_i x_i, 1\leq i\leq n\} \simeq K^{*n} \;\;\;\; ({\rm Lemma} \;
\ref{t19Sep06}).$$ Therefore, 
\begin{equation}\label{SGhTn}
\CS = G\cdot (h_1, \ldots , h_n) \simeq G/ \mathbb{T}^n.
\end{equation}
Note that $\St (h_1, \ldots , h_n) = \{ \s \in G \, | \, \s h_i=
h_i\s,1\leq i \leq n\}$.
\begin{lemma}\label{t19Sep06}
Let $K$ be a reduced commutative ring with $\frac{1}{2}\in K$.
Then $\St (h_1, \ldots , h_n) = \mathbb{T}^n$.
\end{lemma}

{\it Proof}. Clearly, the torus is a subgroup of the stabilizer.
We have to show that each element $\s$ of the stabilizer belongs
to the torus. Since the automorphism $\s$ commutes with all the
$h_i$, the automorphism $\s$ respects  eigenspaces of $h_i$ (i.e.
$\ker (h_i)$ and $\ker (h_i-1)$) and their intersections. In
particular, for each $ i=1, \ldots , n$, the vector space
\begin{equation}\label{CKxi}
\ker (h_1) \cap \cdots \cap \ker (h_{i-1}) \cap \ker (h_i-1)\cap
\ker (h_{i+1}) \cap \cdots \cap \ker (h_n) = Kx_i
\end{equation}
is $\s$-invariant, i.e. $\s (x_i)= \l_ix_i$ for some $\l_i\in
K^*$, and so $\s \in \mathbb{T}^n$.  $\Box $

Let $\CA$ be the set of all the $n$-tuples $(x_1', \ldots , x_n')$
of canonical generators for the $K$-algebra $\L_n(K)$ ($x_i'^2=0$
and $x_i'x_j'=-x_j'x_i'$). Clearly, all $x_i'\in \gm$ since $K$ is
reduced. The group $G$ acts on the set $\CA$ in the obvious way:
$\s (x_1', \ldots , x_n')= (\s (x_1'), \ldots , \s (x_n'))$. The
action is transitive and the stabilizer of each point is trivial
(by the very definition of the group $G$).

\begin{lemma}\label{e19Sep06}
Let $K$ be a reduced commutative ring with $\frac{1}{2}\in K$,
$(s_1, \ldots , s_n)\in \CS$, and $\CK_i = Kx_i'$, $1\leq i \leq
n$. Then the element $(x_1', \ldots , x_n')$ belongs to the set
$\CA$, and $s_1=x_1'\frac{\der}{\der x_1'}, \ldots
,s_n=x_1'\frac{\der}{\der x_n'}$.
\end{lemma}

{\it Proof}.   By  Proposition \ref{m19Sep06}, $\d (\gm )
\subseteq \gm$ for any $K$-derivation $\d$ of $\L_n$. In
particular, $s_1(\gm ) \subseteq \gm , \ldots ,s_n(\gm ) \subseteq
\gm$. Using $n-1$ times (\ref{K1i}), we have
\begin{eqnarray*}
 \L_n&=& K_1+x_1'K_1=( K_{1,2}+ x_2'K_{1,2})+ x_1'( K_{1,2}+ x_2'K_{1,2})\\
 &=&K_{1,2}+x_2'K_{1,2}+x_1'K_{1,2}+x_1'x_2'K_{1,2}=\cdots \\
&=& \sum_{\alpha \in \CB_n} x'^\alpha K_{1, \ldots , n}  =
\sum_{\alpha \in \CB_n} x'^\alpha K = \sum_{\alpha \in \CB_n}
Kx'^\alpha
\end{eqnarray*}
since $ K_{1, \ldots , n}=K$. Since $\L_n(K)$ is a free module of
rank $2^n$ over the commutative ring $K$ with identity,
$x'^\alpha \neq 0$ for all $\alpha$, and the sums above are the
direct sums, i.e. $\L_n= \oplus_{\alpha \in \CB_n} Kx'^\alpha$.

For each $i$ and $\alpha$, $s_i( x'^\alpha)= \alpha_ix'^\alpha$
and $ s_i(y_{\alpha , \nu })=\alpha_iy_{\alpha , \nu }$ where
$y_{\alpha , \nu }:= x_{\nu (1)}'^{\alpha_{\nu (1)}}\cdots x_{\nu
(n)}'^{\alpha_{\nu (n)}}$ and $\nu \in S_n$ ($S_n$ is the
symmetric group). Therefore, $Kx'^\alpha := Kx_1'^{\alpha_1}\cdots
x_n'^{\alpha_n}= K x_{\nu (1)}'^{\alpha_{\nu (1)}}\cdots x_{\nu
(n)}'^{\alpha_{\nu (n)}}$ for any permutation $\nu \in S_n$ (since
the sums above are direct). In particular, for each $i\neq j$,
$x_i'x_j'= \l x_j'x_i'$ for some $\l = \l_{ij}\in K$. We claim
that $\l =-1$. For, note that $\gm = (x_1', \ldots , x_n')$ and so
the set $\bx_1':= x_1'+\gm^2, \ldots , \bx_n':= x_n'+\gm^2$ is a
basis for the vector space $\gm / \gm^2$ over $K$. In $\gm^2/
\gm^3$, on the one hand, $\bx_i'\bx_j' = - \bx_j' \bx_i'\neq 0$;
on the other, by taking the equation $x_i'x_j'= \l x_j'x_i'$
modulo $\gm^3$, we have $\bx_i'\bx_j' = \l \bx_j' \bx_i'$; hence
$\l=-1$, as required.

For each $i$, $x_i'^2\in \ker (s_i)$ (since $s_i(x_i'^2) =
2x_i'^2$  and 2 is not an eigenvalue for the idempotent derivation
$s_i$ as $\frac{1}{2}\in K$), hence
  $x_i'^2\in \gm \cap \cap_{i=1}^n \ker (s_i)= \gm
\cap K=0$, i.e. $x_i'^2=0$. This proves that the elements $x_1',
\ldots x_n'$ are {\em canonical} generators for the algebra
$\L_n$. Now, it is obvious that $s_1=x_1'\frac{\der}{\der x_1'},
\ldots ,s_n=x_n'\frac{\der}{\der x_n'}$.  $\Box $

If a group $\CG$ acts on a set $X$ we say that $X$ is a $\CG$-{\em
set}. Let $Y$ be a $\CG$-set. A map $f:X\ra Y$ is called a
$\CG$-map if $f(gx) = gf(x)$  for all $x\in X$ and $g\in \CG$. A
$\CG$-{\em isomorphism} is a $\CG$-map which is a bijection. The
torus $\mathbb{T}^n$ acts on the set $\CA$ by the rule $(\l_i)
(x_i') := (\l_ix_i')$. Let $ \CA /\mathbb{T}^n$ be the set of all
$\mathbb{T}^n$-orbits.

\begin{corollary}\label{a19Sep06}
Let $K$ be a reduced commutative ring with $\frac{1}{2}\in K$.
Then
\begin{enumerate}
\item The map $\CS \ra \CA /\mathbb{T}^n$, $ (s_1, \ldots, s_n)
\mapsto \mathbb{T}^n (x'_1, \ldots , x'_n)$, is a $G$-isomorphism
with the inverse $\mathbb{T}^n(x_1', \ldots , x_n')\mapsto
(x_1'\frac{\der}{\der x_1'}, \ldots ,x_n'\frac{\der}{\der x_n'})$.
\item In particular, $G$ acts transitively on the set $\CS$.
\end{enumerate}
\end{corollary}

{\it Proof}. 1. This follows directly from Lemma \ref{e19Sep06}.

2. The group $G$ acts transitively on the set $\CA$, hence it does
on the set $\CS$, by statement 1. $\Box $

For $n=1$, Corollary \ref{a19Sep06} gives $\CS = \{ x_1\der_1\}$
since $G=\mathbb{T}$.

Let $\CG$ be a group and $M$ be a $\CG$-module. We say that $M$ is
a {\em faithful} $\CG$-module (or the group $\CG$ acts {\em
faithfully} on $M$) if the map $\CG \ra {\rm End} (M)$ is
injective.

\begin{theorem}\label{26Sep06}
Let $K$ be a reduced commutative ring with $\frac{1}{2}\in K$,
$n\geq 1$. Then
\begin{enumerate}
\item  $\Der_K(\L_n)$ is a faithful $G$-module iff $n\geq 2$.\item
The $G$-module $\Der_K(\L_n)$ is not simple.
\end{enumerate}
\end{theorem}

{\it Proof}. 1. For $n=1$, $\Der_K(\L_1) = Kx_1\der_1$ (Corollary
\ref{cc13Sep06}) and $G=\mathbb{T}={\rm St} (x_1\der_1)$.
Therefore, $\Der_K(\L_1)$ is not a faithful $G$-module. So, let
$n\geq 2$.
  Suppose that an element $\s \in G$ acts trivially on
$\Der_K(\L_n)$, i.e. $\s \d \s^{-1}=\d$ for all $\d \in
\Der_K(\L_n)$. We have to show that $\s = e$, the identity element
of $G$. By Lemma \ref{t19Sep06}, $\s = \s_\l \in \mathbb{T}^n$ for
some $\l  \in K^{*n}$. For each $i=1, \ldots , n$, $\ad ( x_i) =
\s \ad (x_i)\s^{-1} = \ad (\s (x_i))= \l_i\ad (x_i)$, hence
$\l_i=1$ (choose $j$ such that $j\neq i$; then $0=(\l_i-1) \ad
(x_i)(x_j) = 2(\l_i-1) x_ix_j$, and so $\l_i-1=0$), i.e. $\s =e$.

2. If $n\geq 1$, then the $G$-module $\Der_K(\L_n)= D_0 \oplus
D_+$
 contains the  proper
submodule $D_+$, and so $\Der_K(\L_n)$ is not a simple $G$-module.
 $\Box $


\section{Skew derivations of the Grassmann
rings}\label{SKDGA}

Let $K$ be a commutative ring. Recall that the Grassmann
$K$-algebra $\L_n = \Lnz \oplus \Ln1$ is a $\Z_2$-graded algebra,
$ \Lnz = \Lnev$ and $\Ln1 = \Lnod$. Each element $a$ of $\L_n$ is
a unique sum $a= \az + \a1$ with $\az\in \Lnz$ and $\a1\in \Ln1$.
We also use the alternative notation: $a= a^{ev}+a^{od}$ where
$a^{ev}:=\az$ and $a^{od}:= \a1 $.

Recall that a $K$-linear map $\d : \L_n\ra \L_n$ is called a {\em
(left)  skew $K$-derivation}, if for any $b_s\in \L_{n,s}$ and $
b_t\in \L_{n,t}$ (where $s,t\in \Z_2$),
$$ \d ( b_sb_t) = \d (b_s) b_t+ (-1)^s b_s\d (b_t).$$
The set of all skew derivations $\SDer_K(\L_n)$ is a left
$Z(\L_n)$-module and a left $\Lnev$-module since $\Lnev \subseteq
Z(\L_n)$.

Consider the sets of even and odd skew $K$-derivations of
$\L_n(K)$:
\begin{eqnarray*}
\SDerev & :=& \{ \d \in \SDer_K(\L_n )\, | \, \d (\Lni ) \subseteq \Lni , \; \overline{i} \in \Z_2 \},\\
 \SDerod & :=& \{ \d \in \SDer_K(\L_n )\, | \, \d (\Lni ) \subseteq \L_{n, \overline{i}+\overline{1}} , \; \overline{i} \in \Z_2
 \}.
\end{eqnarray*}
So, even skew derivations are precisely the skew derivations that
respect $\Z_2$-grading of $\L_n$, and the odd skew derivations are
precisely the skew derivations that reverse it. The set of odd and
even skew derivations are left $\Lnz$-modules. For each element
$a\in \L_n$, one can attach, so-called, the {\em inner} skew
$K$-derivation of $\L_n$: $\sad (a) : b_s\mapsto ab_s-(-1)^sb_sa$,
where $b_s\in \L_{n, s}$, $s\in \Z_2 $.   The set of all inner
skew derivations is denoted by $\ISDer_K(\L_n)$.
 The kernel of the $K$-linear map $\sad : \L_m\ra \ISDer_K(\L_n)$,
  $a\mapsto \sad (a)$,  is equal to $\ker (\sad )= \Lnod
  +Kx_1\cdots x_n$, and
  $$ \L_n/\ker (\sad )= \Lnev \oplus \Lnod/ (\Lnod
  +Kx_1\cdots x_n)\simeq \Lnev / \Lnev \cap Kx_1\cdots x_n\simeq
  \Lnsev.$$
By the Homomorphism Theorem, the map 
\begin{equation}\label{Ln1s}
\Lnsev \ra \ISDer_K(\L_n), \;\; a\mapsto \sad (a),
\end{equation}
 is a bijection and $\ISDer_K(\L_n)= \{ \sad (a) \, | \,
 a\in \Lnsev\}$.

The next theorem describes
 explicitly the sets of all/inner/even and odd skew derivations.
\begin{theorem}\label{a13Sep06}
Suppose that $K$ is a commutative ring with $\frac{1}{2}\in K$.
Then
\begin{enumerate}
\item $\SDer_K(\L_n) = \SDerev \oplus \SDerod$. \item $\SDerod =
\oplus_{i=1}^n \Lnev\der_i$. \item $\SDerev = \ISDer_K(\L_n)= \{
\sad (a) \, | \, \in \Lnsev \}$,  and the map $ \sad : \Lnsev \ra
\ISDer_K(\L_n)$, $a\mapsto \sad (a)$,  is a bijection. \item
$\SDer_K(\L_n)/ \ISDer_K(\L_n)\simeq \SDerod$.
\end{enumerate}
\end{theorem}

{\it Proof}. By (\ref{Ln1s}), the map in statement 3 is a
bijection and $\ISDer_K(\L_n)= \{ \sad (a) \, | \, a\in \Lnsev\}$.
 Since $\SDerev \cap \SDerod = 0$, one has the inclusion
\begin{equation}\label{s1ri}
\SDer_K(\L_n)\supseteq \SDerev \oplus \SDerod .
\end{equation}
Clearly, 
\begin{equation}\label{s2ri}
\SDerod \supseteq \sum_{i=1}^n \Lnev \der_i =\bigoplus_{i=1}^n
\Lnev \der_i .
\end{equation}
 For each $a\in \Lnsev$, $\sad (a) \in
\SDerev$, hence 
\begin{equation}\label{s3ri}
\SDerev \supseteq \ISDer_K(\L_n).
\end{equation}
Note that statement 4 follows from statements 1 and 3. Now, it is
obvious that in order to finish the proof of the theorem it
suffices to prove the next claim.

{\it Claim}. $\SDer_K(\L_n) \subseteq \sum_{i=1}^n \Lnev \der_i
+\ISDer_K(\L_n)$.

Indeed, suppose that the inclusion of the claim holds then
$$ \SDer_K(\L_n) \subseteq \sum_{i=1}^n \Lnev \der_i
+\ISDer_K(\L_n)\subseteq \SDerev \oplus \SDerod,$$ hence statement
1 is true by (\ref{s1ri}). Statement 1 together with inclusions
(\ref{s2ri}) and (\ref{s3ri}) implies statements 2 and 3.

{\it Proof of the Claim}. Let $\d$ be a skew $K$-derivation of
$\L_n$. We have to represent the derivation $\d $ as a sum
$$ \d = \sum_{i=1}^n a_i\der_i + \sad (a), \;\; a_i\in \Lnev , \;\;
a\in \Lnsod.$$ The proof of the claim is constructive.
  By (\ref{u=ueo}), for each $i=1, \ldots , n$, let 
\begin{equation}\label{1Lnz}
u_i:= \d (x_i)= \uev_i + \uod_i, \;\;\; {\rm where} \;\; \uev_i\in
\Lnev, \;\; \uod_i\in \Lnod .
\end{equation}
Let $\der := \sum_{i=1}^n \uev_i\der_i$ and  $ \d':= \d - \der$.
Then $\d' (x_i)= \uod_i$.  Note that $\der \in \sum_{i=1}^n \Lnev
\der_i$. Then  changing $\d$ for $\d'$, if necessary, one may
assume that  all $u_i$ belong to $\Lnod$. Now,  it suffices to
show that $\d = \sad (a)$ for some $a\in \Lnev$. We construct such
an $a$ in several steps.

{\it Step 1}. Let us prove that,  {\em for each $i=1, \ldots , n$,
$u_i= v_ix_i$ for some element } $v_i\in K\lfloor x_1, \ldots ,
\widehat{x_i}, \ldots , x_n\rfloor^{ev}$. Note that $0= \d (0)= \d
(x_i^2) = u_ix_i - x_iu_i=2u_ix_i$ since $u_i$ is odd, and so
$u_ix_i=0$ since $\frac{1}{2}\in K$. It follows that   $u_i=
v_ix_i$ for some element $v_i\in K\lfloor x_1, \ldots ,
\widehat{x_i}, \ldots , x_n\rfloor^{ev}$. If $n=1$ then $v_i\in
K$, and so $\d = \sad (\frac{1}{2}v_i)$, and we are done. So, let
$n\geq 2$.

{\it Step 2}. We claim that, {\em for each  pair} $i\neq j$,
\begin{equation}\label{sast2}
v_i|_{x_j=0}=v_j|_{x_i=0}.
\end{equation}
Evaluating the skew derivation $\d$ at the element $0=
x_ix_j+x_jx_i$ and taking into account that all the elements $u_i$
are odd we get
\begin{eqnarray*}
 0&=&  u_ix_j-x_iu_j+u_jx_i-x_ju_i= 2(u_ix_j+u_jx_i) \\
 &=& 2(v_ix_ix_j+v_jx_jx_i)= 2(v_i-v_j) x_ix_j.
\end{eqnarray*}
This means that $v_i-v_j\in (x_i, x_j)$ since $\frac{1}{2}\in K$,
or, equivalently, $v_i|_{x_i=0, x_j=0}=v_j|_{x_i=0, x_j=0}$. By
Step 1, this equality can be written as (\ref{sast2}).

{\it Step 3}. Note that $\d (x_1)=v_1x_1= \sad (\frac{1}{2}v_1)
(x_1)$ since  $v_1$ is even, and so  $(\d -\sad
(\frac{1}{2}v_1))(x_1)=0$. Changing $\d$ for $\d -\sad
(\frac{1}{2}v_1)$ one can assume that $\d (x_1)=0$, i.e. $v_1=0$.
Then, by (\ref{sast2}), $v_i|_{x_1=0}=0$ for all $i=2, \ldots ,
n$, and so $v_i\in (x_1)$ for all $i=2, \ldots , n$.
 The idea of the proof is to continue in this way killing the
 elements $v_i$.  Namely, we are going to prove by induction on
 $k$ that, for each $k$ such that $1\leq k\leq n$, by adding to
 $\d$ a well chosen inner skew derivation we have
\begin{equation}\label{skind}
\d (x_1)=\cdots = \d (x_k)=0, \;\; \d (x_i) \in (x_1 \cdots
x_kx_i), \;\; k<i\leq n .
\end{equation}
The case $k=1$ has just been established. Suppose that
(\ref{skind}) holds for $k$, we have to  prove the same statement
for $k+1$. Note that $v_{k+1} x_{k+1} = \d ( x_{k+1})\in
(x_1\cdots x_kx_{k+1})$  (see Steps 1 and 2), hence $v_{k+1} =
x_1\cdots x_kv$ for some $v\in K\lfloor x_{k+2} , \ldots ,
x_n\rfloor$. Consider the derivation $\d':= \d -\sad
(\frac{1}{2}v_{k+1})$. For each $i=1, \ldots , k$, $\d'(x_i)= \d
(x_i) =0$ as $v_{k+1} \in ( x_1\cdots x_k)$; and
\begin{eqnarray*}
\d'( x_{k+1}) &=& v_{k+1} x_{k+1}- \frac{1}{2}(v_{k+1}
x_{k+1}-(-1)x_{k+1}v_{k+1})\\
&=&v_{k+1} x_{k+1}- \frac{1}{2}(v_{k+1} x_{k+1}+v_{k+1}x_{k+1})=0.
\end{eqnarray*}
So, we have proved  the first part of (\ref{skind}) for $k+1$,
namely, that
$$ \d' (x_1)=\cdots = \d' (x_{k+1})=0.$$
So, changing $\d$ for $\d'$ one can assume that
$$ \d (x_1)=\cdots = \d (x_{k+1})=0.$$
These conditions imply that $v_1=\cdots = v_{k+1}=0$ (by Step 1).
If $n=k+1$, we are done. So, let $k+1<n$.  Then, by (\ref{sast2}),
for each $i>k+1$, $v_i \in \cap_{j=1}^{k+1} (x_j)= (x_1\cdots
x_{k+1})$, hence $\d (x_i) = v_ix_i \in ( x_1\cdots x_{k+1} x_i)$,
and we are done.
 By induction, (\ref{skind})
is true for all $k$. In particular, for $k=n$ one has $\d =0$.
This means that $\d$ is an inner skew  derivation, as required.
$\Box $

By Theorem \ref{a13Sep06}, any skew $K$-derivation $\d$ of $\L_n$
is a unique sum $\d = \dev +\dod$ of {\em even} and {\em odd}
 skew derivations, and $\dev := \frac{1}{2}\sad (a)$ for a {\em unique}
element $a\in \Lnsev$. The next corollary describes explicitly the
skew derivations $\dev$ and $\dod$.

\begin{corollary}\label{ca13Sep06}
Let $K$ be a commutative ring with $\frac{1}{2}\in K$, $\d$ be a
skew  $K$-derivation of $\L_n(K)$, and, for each $i=1, \ldots ,
n$,  $\d (x_i) = \uev_i + \uod_i$ for unique  elements $\uev_i\in
\Lnev $ and $\uod_i\in \Lnod $. Then
\begin{enumerate}
\item $\dod =\sum_{i=1}^n \uev_i\der_i$, and \item $\dev =
\frac{1}{2}\sad (a)$ where the {\em unique} element $a\in \Lnsev$
is given by the formula
$$ a= \sum_{i=1}^{n-1} x_1\cdots x_i\der_i\cdots \der_1\der_{i+1}
(\uod_{i+1})+\der_1 (\uod_1).$$
\end{enumerate}
\end{corollary}

{\it Proof}. 1. This statement has been proved already  in the
proof of Theorem \ref{a13Sep06}.

2. For each $i=1, \ldots , n$, on the one hand $\dev (x_i) = (\d -
\dod )(x_i) = \uod_i$; on the other, $ \dev (x_i) =
\frac{1}{2}\sad (a) (x_i)=\frac{1}{2}2x_ia= x_ia$. So, the element
$a$ is a solution to the system of  equations
$$\begin{cases}
x_1a=\uod_1 \\
x_2a=\uod_2 \\
\;\;\;\; \;\;\;\vdots \\
x_na=\uod_n.
\end{cases}
$$
By Theorem \ref{s14Sep06} and the fact that $a\in \Lnsev$ (i.e.
$a_n=0$), we have
$$ a= \sum_{i=1}^{n-1} x_1\cdots x_i\der_i\cdots \der_1\der_{i+1}
(\uod_{i+1})+\der_1 (\uod_1).\;\;\; \Box$$

{\it Definition}. An ideal $\ga$ of $\L_n$ is called a {\em skew
differential} ideal if $\d (\ga ) \subseteq \ga$ for all $\d \in
\SDer_K(\L_n)$.

\begin{lemma}\label{16Mar7}
Let $K$ be a commutative ring with $\frac{1}{2}\in K$, $\CI (K)$
be the set of ideals of the ring $K$, $\SDI (\L_n)$ be the set of
all skew differential ideals of $\L_n$. Then the map
$$ \CI (K) \ra \SDI (\L_n), \;\; I\mapsto I\L_n,$$
is a bijection. In particular, if $K$ is a field of characteristic
$\neq 2$, then $\L_n$ is a skew differentiably simple algebra,
i.e. $0$ and $\L_n$ are the only skew differential ideals of
$\L_n$.
\end{lemma}

{\it Proof}. The map $I\mapsto I\L_n$ is well-defined an
injective. It remains to prove that it is surjective. Let $\ga$ be
a skew differential ideal of $\L_n$. First, let us show that
$$\ga = \oplus_{\alpha \in \CB_n} (\ga \cap Kx^\alpha ) ,$$ i.e. if $a=
\sum_{\alpha \in \CB_n} \l_\alpha x^\alpha\in \ga$, $\l_\alpha\in
K$, then all $\l_\alpha x^\alpha\in \ga$. The case $a=0$ is
trivial. So, let $a\neq 0$ and $i:= \max \{ |\alpha | \, | \,
\l_\alpha \neq 0\}$. We use induction on $i$. The case $i=0$ is
obvious. So, let $i>0$. Then, $\l_\alpha = \der^\alpha (a)\in \ga$
for each $\alpha $ such that $|\alpha | = i$ where $\der^\alpha :=
\der_n^{\alpha_n} \der_{n-1}^{\alpha_{n-1}}\cdots
\der_1^{\alpha_1}$. Applying induction to the element
$a-\sum_{|\alpha | =i}\l_\alpha x^\alpha\in \ga$, we get the
result. So, $\ga = \oplus_{\alpha \in \CB_n} \ga_\alpha x^\alpha$
for some ideals $\ga_\alpha $ of $K$. Let $I:= \ga_0$. On the one
hand, $I\L_n\subseteq \ga$, and so $I\subseteq \ga_\alpha$ for all
$\alpha\in \CB_n$. On the other, $\ga_\alpha =
\der^\alpha(\ga_\alpha  x^\alpha ) \subseteq I$, hence $\ga =
I\L_n$. So, the map $I\mapsto I\L_n$  is a surjection.  $\Box $


\section{Normal elements of the Grassmann algebras}\label{NEGA}

In this section, it is proved that the set of `generic' normal
non-units forms no more that two orbits under the action of the
group $G:= {\rm Aut}_K(\L_n )$ (Theorem \ref{5Nov06}). The
stabilizers of elements from each orbit are found (Lemma
\ref{c5Nov06} and Lemma \ref{s10Nov06}).

In this section, $K$ is a {\em reduced commutative} ring with
$\frac{1}{2}\in K$ and $n\geq 2$. Recall than an element $r$ of a
ring $R$ is called a {\em normal} element if $rR=Rr$. Each unit is
a normal element.

Recall that $G:= \Aut_K(\L_n)$ is the group of $K$-automorphisms
of $\L_n$. Consider some of its subgroups:

\begin{itemize}
 \item $\O := \{ \o_{1+a}\, | \, a\in \L_n'^{od} \}$, where $\o_u :
 \L_n\ra \L_n$, $ x\mapsto uxu^{-1}$, is an inner automorphism.
 \item $\G := \{ \g_b \, | \, \g_b(x_i) = x_i+b_i,
\; b_i \in \Lnod \cap \gm^3, i=1, \ldots , n\}$, $ b=(b_1, \ldots
, b_n)$,  \item $\GL_n(K)^{op}:= \{ \s_A\, | \, \s_A(x_i)=
\sum_{j=1}^n a_{ij}x_j, \; A=(a_{ij})\in \GL_n(K)\}$.
\end{itemize}
For each $a\in \Lnod$ and $x\in \L_n$, $\o_{1+a} (x) = x+[a,x]$
(Lemma 2.8.(3), \cite{jacgras}). Note that
 $G= (\O \rtimes \G ) \rtimes \GL_n(K)^{op}$ (Theorem 2.14, \cite{jacgras}). So, each element $\s
\in G$ has the unique presentation as the product $\s = \o_{1+a}
\g_b\s_A $ where $\o_{1+a} \in \O$ ($a\in \L_n'^{od}$), $\g_b \in
\G$, $\s_A\in \GL_n(K)^{op}$ where $\L_n'^{od} := \oplus_i
\L_{n,i}$ and  $i$ runs through all odd natural numbers such that
$1\leq i\leq n-1$. For more information on the group $G$, the
reader is refereed to \cite{Berezin67} and \cite{Djokovic78} where
$K$ is a field, and to \cite{jacgras} where $K$ is a commutative
ring.

\begin{theorem}\label{M30Sep06}
\cite{jacgras} Let  $K$ be a  reduced commutative ring with
$\frac{1}{2}\in K$. Then each element $\s \in G$ is a unique
product $\s = \o_{1+a} \g_b\s_A$ where $a\in \L_n'^{od}$ and
\begin{enumerate}
\item $\s (x) = Ax +\cdots $ (i.e. $\s (x) \equiv Ax\mod \gm$) for
some $A\in \GL_n(K)$, \item $b= A^{-1} \s(x)^{od}-x$, and \item
$a= -\frac{1}{2}\g_b (\sum_{i=1}^{n-1} x_1\cdots x_i\der_i\cdots
\der_1\der_{i+1} (a_{i+1}') +\der_1(a_1'))$ where $a_i':= (A^{-1}
\g_b^{-1} (\s (x)^{ev}))_i$, the $i$'th component of the
column-vector $A^{-1} \g_b^{-1} (\s (x)^{ev})$.
\end{enumerate}
\end{theorem}

{\it Remark}. In the above theorem the following abbreviations are
used

 $x=\begin{pmatrix}
  x_1\\
 \vdots  \\
x_n \\
\end{pmatrix}$, $b=\begin{pmatrix}
  b_1\\
 \vdots  \\
b_n \\
\end{pmatrix}$, $\s (x)=\begin{pmatrix}
  \s (x_1)\\
 \vdots  \\
\s (x_n) \\
\end{pmatrix}$, $\s (x)^{ev}=\begin{pmatrix}
  \s (x_1)^{ev}\\
 \vdots  \\
\s (x_n)^{ev} \\
\end{pmatrix}$, $\s (x)^{od}=\begin{pmatrix}
  \s (x_1)^{od}\\
 \vdots  \\
\s (x_n)^{od} \\
\end{pmatrix}$, any element $u\in \L_n$ is a unique sum $u=
u^{ev}+u^{od}$ of its even and odd components. Note that the
inversion formula for $\g_b^{-1}$  is given in \cite{jacgras}.

Let $\CN$ be the set of all  normal elements of the Grassmann
algebra $\L_n = \L_n(K)$ and let $\CU$ be the set of all units of
$\L_n$. Then $\CU \subseteq \CN$. By (\ref{smi=mi}), the set $\CN$
is a  disjoint union of its $G$-invariant subsets,
$$\CN = \cup_{i=0}^n\CN_i, \;\; \CN_i:= \{ a\in \CN \, | \, a=
a_i+\cdots , 0\neq a_i\in \L_{n,i}\}.$$ Clearly, $\CN_0= \CU$.
Similarly, by (\ref{smi=mi}), the set $\CU$ is a disjoint union of
its $G$-invariant subsets $\CU_i$,
$$\CU = \cup_{i=0}^n\CU_i, \;\; \CU_i:= \{ a\in \CU \, | \, a=a_0+
a_i+\cdots , a_0\in K^*,  0\neq a_i\in \L_{n,i}\}.$$

\begin{lemma}\label{k5Nov06}
$\Lnev \cup \Lnod \subseteq \CN$.
\end{lemma}

{\it Proof}. It is obvious. $\Box $

The next result shows that `generic' normal non-unit elements of
$\L_n$ (i.e.  the set $\CN_1$) form a single $G$-orbit  if $n$ is
even, and two $G$-orbits if $n$ is odd.

\begin{theorem}\label{5Nov06}
Let $K$ be a field of characteristic $\neq 2$ and $\L_n =
\L_n(K)$. Then\begin{enumerate} \item $\CN_1= Gx_1$ if $n$ is
even.\item $\CN_1= Gx_1\cup G(x_1+x_2\cdots x_n)$ is the disjoint
union of two orbits if $n$ is odd.
\end{enumerate}
\end{theorem}

{\it Proof}. The elements $x_1$ and $y:=x_1+x_2\cdots x_n$ are
normal. First, let us prove that if $n$ is odd then the orbits
$Gx_1$ and $Gy$ are distinct. Suppose that they coincide, i.e.
$y=\s (x_1)$ for some automorphism $\s \in G$, we seek a
contradiction. By Theorem \ref{M30Sep06}, $\s = \o_{1+a}\g_b\s_A$.
By taking the equality $\s (x_1) = y$ modulo the ideal $\gm^2$, we
have $\s_A(x_1)= x_1$, hence $$x_1+x_2\cdots x_n= y= \s (x_1) =
\o_{1+a} \g_b(x_1) = \o_{1+a} (x_1+b_1) = x_1+b_1+[a, x_1+b_1].$$
Equating the odd parts of both ends of the equalities above gives
$x_1= x_1+b_1$, hence $b_1=0$ and $x_2\cdots x_n= [a,x_1]\in
(x_1)$, a contradiction. Therefore, the orbits $Gx_1$ and $Gy$ are
distinct.

It remains to prove that $\CN_1\subseteq Gx_1$ and $\CN_1\subseteq
Gx_1\cup Gy$ in the first and the second case respectively.

Let $a= a_1+a_2+\cdots \in \CN_1$ where all $a_i\in \L_{n,i}$ and
$0\neq a_1\in \L_{n, 1}$. Up to the action of the group
$\GL_n(K)^{op}$, one can assume that $a_1=x_1$. The automorphism
$\g : x_1\mapsto x_1+a_3+a_5+\cdots$, $x_i\mapsto x_i$, $i\geq 2$,
is an element of the group $\G$. Now, $a= \g (x_1) +a^{ev}$ where
$a^{ev} := a_2+a_4+\cdots $ is the even part of the element $a$,
and so $\g^{-1} (a) = x_1+\g^{-1} (a^{ev})$. Note that $\g^{-1}
(a^{ev})$ is an even element of the set $\Lnev \cap \gm^2$.
Therefore, up to the action of the group $\G$, one can assume that
$a= x_1+a^{ev}$ where $a^{ev}$ is an even element of $\gm^2$.
Since $\frac{1}{2}\in K$, the element $a^{ev}$ is the unique sum
$a^{ev}= 2\alpha x_1+\beta$ where $\alpha$ and $\beta$ are
respectively odd and even elements of the Grassmann algebra
$K\lfloor x_2, \ldots , x_n\rfloor$ and $\alpha, \beta \in \gm$.
Applying the inner automorphism $\o_{1-\alpha}=
(\o_{1+\alpha})^{-1}$ to the equality
$$ a= x_1+2\alpha x_1+\beta = x_1+ [ \alpha , x_1]+\beta =
\o_{1+\alpha}(x_1) +\beta = \o_{1+\alpha}(x_1+\beta )$$ we have
the equality $\o_{1-\alpha}(a) = x_1+\beta$. So, up to the action
of the group $\O$, one can assume that $a^{ev} = \beta \in
K\lfloor x_2, \ldots , x_n\rfloor^{ev}_{\geq 2}$. If $\beta =0$
then we are done. So, let $\beta \neq 0$.

{\em Case 1. $\beta x_i=0$ for all} $i=2, \ldots , n$, i.e. $\beta
= \l x_2, \ldots x_n$, $\l \in K$, hence $n$ must be  odd since
$\beta$ is even. In this case, $a= \s_A(y)$ where
$\s_A(x_2):=\l^{-1} x_2$ and $ \s_A(x_i)= x_i$ for all $i\neq 2$,
and we are done.

{\em Case 2. $\beta x_i\neq 0$ for some $i\geq 2$}. We aim to show
that this case is impossible, we seek a contradiction. Let $\beta
= a_{2m} +\cdots$ where $a_{2m}\in K\lfloor x_2, \ldots ,
x_n\rfloor_{2m}$, $a_{2m}x_i\neq 0$,  and the three dots mean
higher terms with respect to the $\Z$-grading of the Grassmann
algebra $\L_n$. The element $a= x_1+a_{2m}+\cdots$ is normal, and
so $ax_i= ba$ for some element $b= b_0+b_1+\cdots \in \L_n$ where
$b_i\in \L_{n,i}$. In more detail,
$$(x_1+a_{2m}+\cdots ) x_i= (b_0+b_1+\cdots )(x_1+a_{2m}+\cdots
).$$ Clearly, $b_0=0$. Equating the homogeneous components of
degrees $1,\ldots , 2m+1$ of both sides of the equality we have
the system of equations:
$$\begin{cases}
b_1x_1= x_1x_i,\\
b_2x_1=0, \\
\;\;\;\; \;\;\;\vdots \\
b_{2m-1}x_1=0,\\
b_1a_{2m} +b_{2m}x_1= a_{2m} x_i.
\end{cases}
$$
The first equation gives $b_1= -x_i+\l x_1$ for some $\l \in K$.
By taking the last equation modulo the ideal $(x_1)$ of $\L_n$ we
have the following equality in the Grassmann algebra $K\lfloor
x_2, \ldots , x_n\rfloor$: $-x_ia_{2m}= a_{2m} x_i$, hence
$2a_{2m} x_i=0$. Dividing by 2, we have the equality
$a_{2m}x_i=0$, which contradicts to the assumption that
$a_{2m}x_i\neq 0$.   $\Box $

Let $\Stab (x_1):= \{ \s \in G\, | \, \s (x_1)=x_1\}$ be the {\em
stabilizer} of the element $x_1$ in $G$ and $\CN_1$ be the set of
normal elements as in Theorem \ref{5Nov06}. By  Theorem
\ref{5Nov06}.(1), the map 
\begin{equation}\label{GStx1}
G/ \Stab (x_1) \ra \CN_1, \;\; \s \Stab (x_1) \mapsto \s (x_1),
\end{equation}
is a bijection where $K$ is a field of characteristic $\neq 2$.

The next lemma describes the stabilizer $\Stab (x_1)$.
\begin{lemma}\label{c5Nov06}
Let $K$ be a
reduced commutative ring with $\frac{1}{2}\in K$, $\O_{x_1}:= \O
\cap \Stab (x_1) = \{ \o_{1+a} \, | \, a\in (x_1)\cap \Lnod \}$,
$\G_{x_1} := \G \cap \Stab (x_1)=\{ \g_b\, | \, b= (0, b_2, \ldots
, b_n) \in \L_{n, \geq 3}^{od} \}$, and $ \GL_n(K)^{op}_{x_1} :=
\GL_n(K)^{op}\cap \Stab (x_1)=\{ \s_A\, | \, \s_A(x_1) = x_1, A\in
\GL_n(K)\}$. Then
$$ \Stab (x_1) = \O_{x_1} \G_{x_1} \GL_n(K)^{op}_{x_1}= (\O_{x_1}\rtimes
\G_{x_1})\rtimes  \GL_n(K)^{op}_{x_1}.$$
\end{lemma}

{\it Proof}. The last equality follows from the equality $G= (\O
\rtimes \G )\rtimes  \GL_n(K)^{op}$ (Theorem 2.14, \cite{jacgras})
 provided the equality before is true. Let $\CM := \O_{x_1}
\G_{x_1} \GL_n(K)^{op}_{x_1}$. Then $\Stab (x_1) \supseteq \CM$.
It remains to prove the reverse inclusion. Let $\s \in \Stab
(x_1)$. By Theorem \ref{M30Sep06}, $\s = \o_{1+a} \g_b\s_A$. Since
$\s (x_1) \equiv \s_A(x_1)\mod \gm^2$ and $\s (x_1) = x_1$ we must
have $\s_A\in \GL_n(K)^{op}_{x_1}$. Now,
$$ x_1= \s (x_1) = \o_{1+a}(x_1+ b_1) = x_1+b_1+[a, x_1+b_1]=
x_1+b_1+ 2a(x_1+b_1),$$and equating the odd parts of the elements
at both ends of the equalities we see that $x_1= x_1+b_1$, hence
$b_1=0$, and so $\g_b\in \G_{x_1}$. Putting $b_1=0$ in the
equalities above gives $x_1= x_1+2ax_1$, hence $0=
x_1^{ev}=(x_1+2ax_1)^{ev} = 2ax_1$, and so $a\in (x_1)$. This
means that $\s \in \O_{x_1}$, as required.  $\Box $

\begin{corollary}\label{d5Nov06}
Let $K$ be a
reduced commutative ring with $\frac{1}{2}\in K$, and $I$ be a
non-empty subset of $\{ 1, \ldots , n\}$. Then
$$ \cap_{i\in I} \Stab (x_i) = (\cap_{i\in I}\O_{x_i})\cdot (\cap_{i\in I} \G_{x_i})\cdot ( \cap_{i\in I}\GL_n(K)^{op}_{x_i})
= (\cap_{i\in I}\O_{x_i})\rtimes (\cap_{i\in I}\G_{x_i})\rtimes (
\cap_{i\in I}\GL_n(K)^{op}_{x_i}).$$
\end{corollary}

{\it Proof}. This follows from Lemma \ref{c5Nov06} (and uniqueness
of the decomposition $ \s = \o_{1+a} \g_b\s_A$). $\Box $

If $K$ is a field of characteristic $\neq 2$ and $n$ is an odd
number then, by Theorem \ref{5Nov06}.(2), the map (where
$y:=x_1+x_2\cdots x_n$)  
\begin{equation}\label{GSty}
G/\Stab (x_1) \cup G/ \Stab (y) \ra \CN , \;\; \s \Stab
(x_1)\mapsto \s (x_1), \;\ \tau \Stab (y)\mapsto \tau (y),
\end{equation}
is a bijection. The next lemma describes the stabilizer $\Stab (
x_1+x_2\cdots x_n)$.

\begin{lemma}\label{s10Nov06}
Let $K$ be a
reduced commutative ring with $\frac{1}{2}\in K$, and  $n\geq 3$
be an  odd number. Then $\Stab (x_1+x_2\cdots x_n)= \{
\o_{1+\frac{1}{2}\der_1\g_b\s_A(x_2\cdots x_n)+x_1c} \g_b \s_A\,\,
|\,\, \g_b\in \G_{x_1},\;\; \s_A\in \GL_n(K)^{op}_{x_1}$, $
(1-x_1\der_1)\g_b \s_A(x_2\cdots x_n) = x_2\cdots x_n, \; c\in
K\lfloor x_2, \ldots , x_n\rfloor^{ev}\}$.
\end{lemma}

{\it Proof}. Let $y:=x_1+x_2\cdots x_n$ and $\s \in \Stab (y)$.
Note that  $\s = \o_{1+a} \g_b\s_A$. Since $x_1\equiv y \equiv \s
(y) \equiv \s (x_1)\mod \gm^2$, we must have $\s_A\in
\GL_n(K)^{op}_{x_1}$. Then

\begin{eqnarray*}
x_1+x_2\cdots x_n&=&y=\s (y) = \o_{1+a}\g_b(x_1+\s_A(x_2\cdots
x_n)) = \o_{1+a} (x_1+b_1+\g_b\s_A(x_2\cdots x_n))\\
& =& x_1+b_1+ \g_b \s_A( x_2\cdots x_n) + [ a, x_1+b_1].
\end{eqnarray*}
Equating the odd parts of the beginning and the end of the series
of equalities above we obtain $x_1= x_1+b_1$, hence $b_1=0$, i.e.
$\g_b\in \G_{x_1}$, and then 
\begin{equation}\label{2x1a}
\g_b\s_A(x_2\cdots x_n) - x_2\cdots x_n = 2x_1a.
\end{equation}
Each element $u\in \L_n$ is a unique sum $u= x_1\alpha + \beta$
for unique elements $\alpha , \beta \in K\lfloor x_2, \ldots ,
x_n\rfloor$. Clearly, $\alpha = \der_1(u)$ and $\beta =
(1-x_1\der_1) (u)$.  The odd element $a$ is a unique sum $a=
x_1c+d$ for some elements $c\in K\lfloor x_2, \ldots ,
x_n\rfloor^{ev}$ and $ d\in K\lfloor x_2, \ldots ,
x_n\rfloor^{od}_{\geq 1}$. The equation (\ref{2x1a}) can be
written as follows $\g_b\s_A (x_2\cdots x_n) - x_2\cdots x_n =
2x_1d$. This equality is equivalent to two equalities $d=
\frac{1}{2}\der_1(\g_b\s_A (x_2\cdots x_n)- x_2\cdots x_n)=
\frac{1}{2}\der_1\g_b\s_A(x_2\cdots x_n)$ and $0=(1-x_1\der_1)
(\g_b\s_A (x_2\cdots x_n)- x_2\cdots x_n)= (1-x_1\der_1)\g_b\s_A
(x_2\cdots x_n)- x_2\cdots x_n$. This finishes the proof of the
lemma. $\Box$

Let  $\CU_1'$ be  the image of the injection $K^*\times \CN_1\ra
\CU_1$, $ (\l , u ) \mapsto \l (1+u)$. The next corollary follows
from Theorem \ref{5Nov06}.

\begin{corollary}\label{u11Nov06}
Let $K$ be a field of characteristic $\neq 2$, and $\L_n =
\L_n(K)$. Then
\begin{enumerate}
\item $\CU_1'= \cup_{\l \in K^*} G\cdot \l (1+x_1)$ is the
disjoint union of orbits if $n$ is even. \item $\CU_1'= \cup_{\l
\in K^*} (G\cdot \l (1+x_1)\cup G\cdot \l (1+x_1+x_2\cdots x_n))$
is the disjoint union of orbits if $n$ is odd.
\end{enumerate}
\end{corollary}


Department of Pure Mathematics

University of Sheffield

Hicks Building

Sheffield S3 7RH

UK

email: v.bavula@sheffield.ac.uk

\end{document}